\documentclass[12pt]{amsart}
\usepackage{amssymb}
\usepackage{amscd,amssymb,amsthm}
\usepackage{graphicx}

\newtheorem{lemma}{Lemma}[section]
\newtheorem{theorem}{Theorem}[section]
\newtheorem{remark}{Remark}[section]
\newtheorem{corollary}{Corollary}[section]

\numberwithin{equation}{section}

\email{taneja@mtm.ufsc.br} \urladdr{http://www.mtm.ufsc.br/$\sim
$taneja}



\begin{document}
\title[Mean Inequalities]{REFINEMENT OF INEQUALITIES AMONG MEANS}

\author{Inder Jeet Taneja}
\address{Inder Jeet Taneja\\
Departamento de Matem\'{a}tica\\
Universidade Federal de Santa Catarina\\
88.040-900 Florian\'{o}polis, SC, Brazil}

\begin{abstract}
In this paper we shall consider some famous \textit{means} such as
\textit{arithmetic}, \textit{harmonic}, \textit{geometric},
\textit{root-square means}, etc. Some new \textit{means} recently
studied are also presented. Different kinds of refinement of
inequalities among these \textit{means} are given.
\end{abstract}

\maketitle

\section{Mean of Order \textit{t}}

Let us consider the following well known \textit{mean of order t}:
\begin{equation}
\label{eq1} B_t (a,b) = \begin{cases}
 {\left( {\frac{a^t + b^t}{2}} \right)^{1 / t},} & {t \ne 0} \\
 {\sqrt {ab} ,} & {t = 0} \\
 {\max \{a,b\},} & {t = \infty } \\
 {\min \{a,b\},} & {t = - \infty } \\
\end{cases}
\end{equation}

\noindent for all $a,b,t \in \mathbb{R},\mbox{ }a,b > 0$.

In particular, we have
\begin{align}
B_{ - 1} (a,b) & = H(a,b) = \frac{2ab}{a + b},\notag\\
B_0 (a,b) & = G(a,b) = \sqrt {ab} ,\notag\\
B_{1 / 2} (a,b) &= N_1 (a,b) = \left( {\frac{\sqrt a + \sqrt b }{2}}
\right)^2,\notag\\
B_1 (a,b) & = A(a,b) = \frac{a + b}{2},\notag\\
\intertext{and} B_2 (a,b) & = S(a,b) = \sqrt {\frac{a^2 +
b^2}{2}}.\notag
\end{align}

The means, $H(a,b)$, $G(a,b)$, $A(a,b)$ and $S(a,b)$ are known in
the literature as \textit{harmonic}, \textit{geometric},
\textit{arithmetic} and \textit{root-square means} respectively. For
simplicity we can call the measure, $N_1 (a,b)$ as
\textit{square-root} mean. It is well know that \cite{beb} the
\textit{mean of order }$s$ given in (\ref{eq1}) is monotonically
increasing in $s$, then we can write
\begin{equation}
\label{eq2}
H(a,b) \leqslant G(a,b) \leqslant N_1 (a,b) \leqslant A(a,b) \leqslant
S(a,b).
\end{equation}

Dragomir and Pearce \cite{drp} (page 242) proved the following
inequality:
\begin{equation}
\label{eq3}
\frac{a^r + b^r}{2} \leqslant \frac{b^{r + 1} - a^{r + 1}}{(r + 1)(b - a)}
\leqslant \left( {\frac{a + b}{2}} \right)^r,
\end{equation}

\noindent for all $a,b > 0$, $a \ne b$, $r \in (0,1)$. In particular
take $r = \frac{1}{2}$ in (\ref{eq3}), we get
\begin{equation}
\label{eq4} \frac{\sqrt a + \sqrt b }{2} \leqslant \frac{2(b^{3 / 2}
- a^{3 / 2})}{3(b - a)} \leqslant \sqrt {\frac{a + b}{2}} , \,\, a
\ne b.
\end{equation}

After necessary calculations in (\ref{eq5}), we get
\begin{equation}
\label{eq5}
\left( {\frac{\sqrt a + \sqrt b }{2}} \right)^2 \leqslant \frac{a + \sqrt
{ab} + b}{3} \leqslant \left( {\frac{\sqrt a + \sqrt b }{2}} \right)\left(
{\sqrt {\frac{a + b}{2}} } \right).
\end{equation}

On the other side we can easily check that
\begin{equation}
\label{eq6}
\left( {\frac{\sqrt a + \sqrt b }{2}} \right)\left( {\sqrt {\frac{a + b}{2}}
} \right) \leqslant \frac{a + b}{2}.
\end{equation}

Finally, the expressions (\ref{eq2}), (\ref{eq5}) and (\ref{eq6})
lead us to the following inequality:
\begin{equation}
\label{eq7}
H(a,b) \leqslant G(a,b) \leqslant N_1 (a,b)
 \leqslant N_3 (a,b) \leqslant N_2 (a,b) \leqslant A(a,b) \leqslant S(a,b),
\end{equation}

\noindent where
\[
N_2 (a,b) = \left( {\frac{\sqrt a + \sqrt b }{2}} \right)\left( {\sqrt
{\frac{a + b}{2}} } \right),
\]

\noindent and
\[
N_3 (a,b) = \frac{a + \sqrt {ab} + b}{3}.
\]

Moreover, we can write
\[
N_1 (a,b) = \frac{A(a,b) + G(a,b)}{2},
\]
\[
N_2 (a,b) = \sqrt {N_1 (a,b)A(a,b)} ,
\]

\noindent and
\[
N_3 (a,b) = \frac{2A(a,b) + G(a,b)}{3}.
\]

Thus we have three new means, where $N_1 (a,b)$ appears as a natural
way. The $N_2 (a,b)$ can be seen in Taneja \cite{tan1,tan2} and the
mean $N_3 (a,b)$ is known as Heron's mean \cite{bob}. Some studies
on it can be seen in Zhang and Wu \cite{zhw}.

\section{Difference of Means and Their Convexity}

Let us consider the following \textit{difference of means}:
\begin{equation}
\label{eq8}
M_{SA} (a,b) = S(a,b) - A(a,b),
\end{equation}
\begin{equation}
\label{eq9}
M_{SN_2 } (a,b) = S(a,b) - N_2 (a,b),
\end{equation}
\begin{equation}
\label{eq10} M_{SN_3 } (a,b) = S(a,b) - N_3 (a,b),
\end{equation}
\begin{equation}
\label{eq11}
M_{SN_1 } (a,b) = S(a,b) - N_1 (a,b),
\end{equation}
\begin{equation}
\label{eq12}
M_{SG} (a,b) = S(a,b) - G(a,b),
\end{equation}
\begin{equation}
\label{eq13}
M_{SH} (a,b) = S(a,b) - H(a,b),
\end{equation}
\begin{equation}
\label{eq14}
M_{AN_2 } (a,b) = A(a,b) - N_2 (a,b),
\end{equation}
\begin{equation}
\label{eq15}
M_{AG} (a,b) = A(a,b) - G(a,b),
\end{equation}
\begin{equation}
\label{eq16}
M_{AH} (a,b) = A(a,b) - H(a,b),
\end{equation}
\begin{equation}
\label{eq17}
M_{N_2 N_1 } (a,b) = N_2 (a,b) - N_1 (a,b),
\end{equation}

\noindent and
\begin{equation}
\label{eq18} M_{N_2 G} (a,b) = N_2
(a,b) - G(a,b).
\end{equation}

We easily check that
\begin{align}
\label{eq19}
M_{AG} (a,b)
& = 2\left[ {N_1 (a,b) - G(a,b)} \right]: = 2M_{N_1 G} (a,b)\\
& = 2\left[ {A(a,b) - N_1 (a,b)} \right]: = 2M_{AN_1 } (a,b)\notag\\
& = 3\left[ {A(a,b) - N_3 (a,b)} \right]: = 3M_{AN_3 } (a,b)\notag\\
& = \frac{3}{2}\left[ {N_3 (a,b) - G(a,b)} \right]: =
\frac{3}{2}M_{N_3 G} (a,b)\notag\\
& = 6\left[ {N_3 (a,b) - N_1 (a,b)} \right]: = 6M_{N_3 N_1 }
(a,b).\notag
\end{align}

Now, we shall prove the convexity of the means
(\ref{eq8})-(\ref{eq18}). It is based on the following lemma.

\begin{lemma} \label{lem21} Let $f:I \subset \mathbb{R}_ + \to
\mathbb{R}$ be a convex and differentiable function satisfying $f(1)
= f^\prime (1) = 0$. Consider a function
\begin{equation}
\label{eq20} \phi _f (a,b) = af\left( {\frac{b}{a}} \right), \,\,
a,b > 0,
\end{equation}

\noindent then the function $\phi _f (a,b)$ is convex in
$\mathbb{R}_ + ^2 $, and satisfies the following inequality:
\begin{equation}
\label{eq21}
0 \leqslant \phi _f (a,b) \leqslant \left( {\frac{b - a}{a}} \right)\phi
_{f}' (a,b).
\end{equation}
\end{lemma}

\begin{proof} It is well known that for the convex and
differentiable function $f$, we have the inequality
\begin{equation}
\label{eq22}
{f}'(x)(y - x) \leqslant f(y) - f(x) \leqslant {f}'(y)(y - x),
\end{equation}

\noindent for all $x,y \in \mathbb{R}_ + $.

Take $y = \frac{b}{a}$ and $x = 1$ in (\ref{eq22}) one gets
\[
{f}'(1)\left( {\frac{b}{a} - 1} \right) \leqslant f\left( {\frac{b}{a}}
\right) - f(1) \leqslant {f}'\left( {\frac{b}{a}} \right)\left( {\frac{b}{a}
- 1} \right),
\]

\noindent or equivalently,
\begin{equation}
\label{eq23}
{f}'(1)\left( {b - a} \right) \leqslant af\left( {\frac{b}{a}} \right) -
af(1) \leqslant a{f}'\left( {\frac{b}{a}} \right)\left( {\frac{b - a}{a}}
\right)
\end{equation}

Since $f(1) = {f}'(1) = 0$, then from (\ref{eq23}) we get
(\ref{eq21}).\\

Now we shall show that the function $\phi _f (a,b)$ is jointly
convex in $a$ and $b$. Since the function $f$ is convex, then for
any $(x_1 ,y_1 )$, $(x_2 ,y_2 ) \,\,  \in \mathbb{R}_ + ^2 $, $0 <
\lambda _1 ,\lambda _2 < 1$, $\lambda _1 + \lambda _2 = 1$ we can
write
\[
f\left( {\frac{\lambda _1 x_1 + \lambda _2 x_2 }{\lambda _1 y_1 + \lambda _2
y_2 }} \right) = f\left( {\frac{\lambda _1 y_1 x_1 }{y_1 \left( {\lambda _1
y_1 + \lambda _2 y_2 } \right)} + \frac{\lambda _2 y_2 x_2 }{y_2 \left(
{\lambda _1 y_1 + \lambda _2 y_2 } \right)}} \right).
\]
\begin{equation}
\label{eq24}
 \leqslant \frac{\lambda _1 y_1 }{\lambda _1 y_1 + \lambda _2 y_2 }f\left(
{\frac{x_1 }{y_1 }} \right) + \frac{\lambda _2 y_2 }{\lambda _1 y_1 +
\lambda _2 y_2 }f\left( {\frac{x_2 }{y_2 }} \right).
\end{equation}

Multiply (\ref{eq24}) by $\lambda _1 y_1 + \lambda _2 y_2 $ one gets
\[
(\lambda _1 y_1 + \lambda _2 y_2 )f\left( {\frac{\lambda _1 x_1 + \lambda _2
x_2 }{\lambda _1 y_1 + \lambda _2 y_2 }} \right) \leqslant \lambda _1 y_1
f\left( {\frac{\lambda _1 x_1 }{\lambda _1 y_1 }} \right) + \lambda _2 y_2
f\left( {\frac{\lambda _2 x_2 }{\lambda _2 y_2 }} \right),
\]

\noindent i.e.,
\begin{equation}
\label{eq25}
\phi _f \left( {\lambda _1 x_1 + \lambda _2 x_2 ,\lambda _1 y_1 + \lambda _2
y_2 } \right) \leqslant \lambda _1 \phi _f \left( {x_1 ,y_1 } \right) +
\lambda _2 \phi _f \left( {x_2 ,y_2 } \right),
\end{equation}

\noindent for any $(x_1 ,y_1 ),\mbox{ }(x_2 ,y_2 ) \in \mathbb{R}_ +
^2 $. The expression (\ref{eq25}) completes the required proof.
\end{proof}

Now we shall show that the \textit{difference of means} given by
(\ref{eq8})-(\ref{eq18}) are \textit{convex} in $\mathbb{R}_ + ^2 $
Later in Section 3 we shall apply the convexity of these functions
to establish improvement over the inequality (\ref{eq7}).

\begin{theorem} \label{the21} The \textit{difference of means} given by
(\ref{eq8})-(\ref{eq18}) are \textit{nonnegative} and
\textit{convex} in $\mathbb{R}_ + ^2 $.
\end{theorem}

\begin{proof} We shall write each measure in the form of
generating function according to the measure (\ref{eq20}), and then
give their first and second order derivatives. It is understood that
$x \in (0,\infty )$.\\

\textbf{$\bullet$ For $M_{SA} (a,b)$:}
\begin{align}
f_{SA} (x) &= \sqrt {\frac{x^2 + 1}{2}} - \frac{x + 1}{2},\notag\\
{f}'_{SA} (x) & = \frac{x}{\sqrt 2 \sqrt {x^2 + 1} } -
\frac{1}{2},\notag \intertext{and} {f}''_{SA} (x) &= \frac{2}{(2x^2
+ 2)^{3 / 2}}
> 0.\notag
\end{align}

\textbf{$\bullet$ For $M_{SN_3 } (a,b)$:}
\begin{align}
f_{SN_3 } (x) &= \sqrt {\frac{x^2 + 1}{2}} - \frac{x + \sqrt x +
1}{3},\notag\\
{f}'_{SN_3 } (x) & = \frac{6x^{3 / 2} - \left( {2\sqrt x + 1}
\right)\sqrt {2(x^2 + 1)} }{6\sqrt {2x(x^2 + 1)} },\notag
\intertext{and} {f}''_{SN_3 } (x) &= \frac{24x^{3 / 2} + \left(
{2x^2 + 2} \right)^{3 / 2}}{12x^{3 / 2}(2x^2 + 2)^{3 / 2}}
 > 0.\notag
\end{align}

\textbf{$\bullet$ For $M_{SN_2 } (a,b)$:}
\begin{align}
f_{SN_2 } (x) & = \frac{2\sqrt {x^2 + 1} - \left( {\sqrt x + 1}
\right)\sqrt {x + 1} }{2\sqrt 2 },\notag\\
{f}'_{_{SN_2 } } (x) & = \frac{4x^{3 / 2}\sqrt {x + 1} - \left( {2x
+ \sqrt x + 1} \right)\sqrt {x^2 + 1} }{4\sqrt {2x(x + 1)(x^2 + 1)}
},\notag \intertext{and} {f}''_{SN_2 } (x) & = \frac{\left( {x^{3 /
2} + 1} \right)\left( {x^2 + 1} \right)^{3 / 2} + 8x^{3 / 2}(x +
1)^{3 / 2}}{8\sqrt 2 \left[ {x(x + 1)(x^2 + 1)} \right]^{3 / 2}}
 > 0.\notag
\end{align}

\textbf{$\bullet$ For $M_{SN_1 } (a,b)$:}
\begin{align}
f_{SN_1 } (x) & = \frac{2\sqrt {2(x^2 + 1)} - \left( {\sqrt x + 1}
\right)^2}{4},\notag\\
{f}'_{_{SN_1 } } (x) & = \frac{4x^{3 / 2} - \left( {\sqrt x + 1}
\right)\sqrt {2(x^2 + 1)} }{4\sqrt {2x(x^2 + 1)} },\notag
\intertext{and} {f}''_{SN_1 } (x) & = \frac{16x^{5 / 2} + x(2x^2 +
2)^{3 / 2}}{8x^{5 / 2}(2x^2 + 2)^{3 / 2}}
 > 0.\notag
\end{align}

\textbf{$\bullet$ For $M_{SG} (a,b)$:}
\begin{align}
f_{SG} (x) & = \sqrt {\frac{x^2 + 1}{2}} - \sqrt x ,\notag\\
{f}'_{SG} (x) &= \frac{\sqrt 2 x^{3 / 2} - \sqrt {x^2 + 1} }{2\sqrt
{x(x^2 + 1)} },\notag \intertext{and} {f}''_{SG} (x) & =
\frac{1}{\sqrt 2 (x^2 + 1)^{3 / 2}} + \frac{1}{4x^{3 / 2}}
 > 0.\notag
\end{align}

\textbf{$\bullet$ For $M_{SH} (a,b)$:}
\begin{align}
f_{SH} (x) & = \sqrt {\frac{x^2 + 1}{2}} - \frac{2x}{x + 1},\notag\\
{f}'_{SH} (x) & = \frac{x(x + 1)^2 - 2\sqrt {2(x^2 + 1)} }{(x +
1)^2\sqrt {2(x^2 + 1)} },\notag \intertext{and} {f}''_{SH} (x) & =
\frac{2\left[ {(x + 1)^3 + 2(2x^2 + 2)^{3 / 2}} \right]}{(x +
1)^3(2x^2 + 2)^{3 / 2}}
 > 0.\notag
\end{align}

\textbf{$\bullet$ For $M_{AN_2 } (a,b)$:}
\begin{align}
f_{AN_2 } (x) & = \frac{2(x + 1) - \left( {\sqrt x + 1} \right)\sqrt
{2(x + 1)} }{4},\notag\\
{f}'_{AN_2 } (x) & = \frac{2\sqrt {2x(x + 1)} - \left( {2x + \sqrt x
+ 1} \right)}{4\sqrt {2(x + 1)} },\notag \intertext{and} {f}''_{AN_2
} (x) & = \frac{x^{3 / 2} + 1}{4x^{3 / 2}(2x + 2)^{3 / 2}}
> 0.\notag
\end{align}

\textbf{$\bullet$ For $M_{AG} (a,b)$:}
\begin{align}
f_{AG} (x) & = \frac{1}{2}(\sqrt x - 1)^2, \notag\\
{f}'_{AG} (x) & = \frac{\sqrt x - 1}{2\sqrt x }\notag
\intertext{and} {f}''_{AG} (x) & = \frac{1}{4x^{3 / 2}}
 > 0.\notag
\end{align}

\textbf{$\bullet$ For $M_{AH} (a,b)$:}
\begin{align}
f_{AH} (x) & = \frac{(x - 1)^2}{2(x + 1)},\notag\\
{f}'_{AH} (x) & = \frac{(x - 1)(x + 3)}{2(x + 1)^2}\notag
\intertext{and} {f}''_{AH} (x) & = \frac{4}{(x + 1)^3}
 > 0.\notag
\end{align}

\textbf{$\bullet$ For $M_{N_2 N_1 } (a,b)$:}
\begin{align}
f_{N_2 N_1 } (x) & = \frac{\left( {\sqrt x + 1} \right)\sqrt {2(x +
1)} - \left( {\sqrt x + 1} \right)^2}{4},\notag\\
{f}'_{N_2 N_1 } (x) & = \frac{2x + \sqrt x + 1 - \left( {\sqrt x +
1} \right)\sqrt {2(x + 1)} }{4\sqrt {2x(x + 1)} },\notag
\intertext{and} {f}''_{N_2 N_1 } (x)  & = \frac{(2x + 2)^{3 / 2} -
2(x^{3 / 2} + 1)}{8x^{3 / 2}(2x + 2)^{3 / 2}}.\notag
\end{align}

Since $(x + 1)^{3 / 2} \geqslant x^{3 / 2} + 1$, $\forall x \in
(0,\infty )$ and $2^{3 / 2} \geqslant 2$, then obviously,
${f}''_{N_2 N_1 } (x) \geqslant 0$, $\forall x \in (0,\infty )$.\\

\textbf{$\bullet$ For $M_{N_2 G} (a,b)$:}
\begin{align}
f_{N_2 G}(x) & = \frac{\left( {\sqrt x + 1} \right)\sqrt {2(x + 1)}
- 4x}{4},\notag\\
{f}'_{_{N_2 G} } (x) & = \frac{2x + 1 + \sqrt x - 2\sqrt {2(x + 1)}
}{4\sqrt {2x(x + 1)} },\notag \intertext{and} {f}''_{N_2 G} (x) & =
\frac{(2x + 2)^{3 / 2} - (x^{3 / 2} + 1)}{4x^{3 / 2}(2x + 2)^{3 /
2}}.\notag \end{align}

Since $(x + 1)^{3 / 2} \geqslant x^{3 / 2} + 1$, $\forall x \in
(0,\infty )$ and $2^{3 / 2} \geqslant 1$, then obviously,
${f}''_{N_2 G} (x) \geqslant 0$, $\forall x \in (0,\infty )$.

We see that in all the cases the generating function $f_{( \cdot )}
(1) = {f}'_{( \cdot )} (1) = 0$ and the second derivative is
positive for all $x \in (0,\infty )$. This proves the
\textit{nonegativity} and \textit{convexity} of the means
(\ref{eq8})-(\ref{eq23}) in $\mathbb{R}_ + ^2 $. This completes the
proof of the theorem.
\end{proof}

\begin{remark} The inequality (\ref{eq7}) also present more
nonnegative differences but we have considered only the convex ones.
\end{remark}

\section{Inequality Among Difference of Means}

In view of (\ref{eq7}), the following inequalities are obviously
true:
\begin{align}
\label{eq26} & M_{SA} (a,b) \leqslant M_{SN_2 } (a,b) \leqslant
M_{SN_3 } (a,b) \leqslant M_{SN_1 } (a,b)
 \leqslant M_{SG} (a,b) \leqslant M_{SH} (a,b),\\
\label{eq27}  & M_{AN_2 } (a,b) \leqslant M_{AN_3 } (a,b) \leqslant
M_{AN_1 } (a,b)
 \leqslant M_{AG} (a,b) \leqslant M_{AH} (a,b),\\
\label{eq28} & M_{N_2 N_3 } (a,b) \leqslant M_{N_2 N_1 } (a,b)
\leqslant M_{N_2 G} (a,b) \leqslant M_{N_2 H} (a,b),\\
\label{eq29} & M_{N_3 N_1 } (a,b) \leqslant M_{N_3 G} (a,b)
\leqslant M_{N_3 H} (a,b), \intertext{and} \label{eq30} & M_{N_1 G}
(a,b) \leqslant M_{N_1 H} (a,b),
\end{align}

In view of (\ref{eq7}), (\ref{eq19}) and (\ref{eq30}), we can easily
check that
\begin{equation}
\label{eq31}
A(a,b) + H(a,b) \leqslant N_1 (a,b) + N_3 (a,b) \leqslant N_1 (a,b) + N_2
(a,b).
\end{equation}

In this section we shall improve the inequalities (\ref{eq7}) and
then compare with the inequalities (\ref{eq26})-(\ref{eq30}). This
refinement is based on the following lemma.

\begin{lemma} \label{lem31} Let $f_1 ,f_2 :I \subset \mathbb{R}_ + \to
\mathbb{R}$ be two convex functions satisfying the assumptions:

(i) $f_1 (1) = f_1 ^\prime (1) = 0$, $f_2 (1) = f_2 ^\prime (1) = 0$;

(ii) $f_1 $ and $f_2 $ are twice differentiable in $\mathbb{R}_ + $;

(iii) there exists the real constants $\alpha ,\beta $ such that $0
\leqslant \alpha < \beta $ and
\begin{equation}
\label{eq32} \alpha \leqslant \frac{f_1 ^{\prime \prime }(x)}{f_2
^{\prime \prime }(x)} \leqslant \beta , \,\, f_2 ^{\prime \prime
}(x) > 0,
\end{equation}

\noindent for all $x > 0$ then we have the inequalities:
\begin{equation}
\label{eq33}
\alpha \mbox{ }\phi _{f_2 } (a,b) \leqslant \phi _{f_1 } (a,b) \leqslant
\beta \mbox{ }\phi _{f_2 } (a,b),
\end{equation}

\noindent for all $a,b \in (0,\infty )$.
\end{lemma}

\begin{proof} Let us consider the functions
\[
k(x) = f_1 (x) - \alpha \mbox{ }f_2 (x)
\]

\noindent and
\[
h(x) = \beta \mbox{ }f_2 (x) - f_1 (x),
\]

\noindent where $\alpha $ and $\beta $are as given by (\ref{eq32}).

In view of item (i), we have $k(1) = h(1) = 0$ and ${k}'(1) =
{h}'(1) = 0$. Since the functions $f_1 (x)$ and $f_2 (x)$ are twice
differentiable, then in view of (\ref{eq32}), we have
\begin{equation}
\label{eq34}
{k}''(x) = f_1 ^{\prime \prime }(x) - \alpha \mbox{ }f_2 ^{\prime \prime
}(x)
 = f_2 ^{\prime \prime }(x)\left( {\frac{f_1 ^{\prime \prime }(x)}{f_2
^{\prime \prime }(x)} - \alpha } \right) \geqslant 0 \,\, ,
\end{equation}

\noindent and
\begin{equation}
\label{eq35}
{h}''(x) = \beta \mbox{ }f_2 ^{\prime \prime }(x) - f_1 ^{\prime \prime }(x)
 = f_2 ^{\prime \prime }(x)\left( {M - \frac{f_1 ^{\prime \prime }(x)}{f_2
^{\prime \prime }(x)}} \right) \geqslant 0,
\end{equation}

\noindent for all $x \in (0,\infty )$.

In view of (\ref{eq34}) and (\ref{eq35}), we can say that the
functions $k( \cdot )$ and $h( \cdot )$, are convex on $I \subset
\mathbb{R}_ + $.

According to (\ref{eq21}), we have
\begin{equation}
\label{eq36}
a\mbox{ }k\left( {\frac{b}{a}} \right) = a\left[ {f_1 \left( {\frac{b}{a}}
\right) - \alpha \mbox{ }f_2 \left( {\frac{b}{a}} \right)} \right] = a\mbox{
}f_1 \left( {\frac{b}{a}} \right) - \alpha a\mbox{ }f_2 \left( {\frac{b}{a}}
\right) \geqslant 0,
\end{equation}

\noindent and
\begin{equation}
\label{eq37}
a\mbox{ }h\left( {\frac{b}{a}} \right) = a\left[ {\beta f_2 \left(
{\frac{b}{a}} \right) - f_1 \left( {\frac{b}{a}} \right)} \right] = \beta
a\mbox{ }f_2 \left( {\frac{b}{a}} \right) - f_1 \left( {\frac{b}{a}} \right)
\geqslant 0,
\end{equation}

Combining (\ref{eq36}) and (\ref{eq37}) we have the proof of
(\ref{eq33}).
\end{proof}

\begin{theorem} \label{the31} The following inequalities among the
\textit{mean differences} hold:
\begin{equation}
\label{eq38}
M_{SA} (a,b) \leqslant \frac{1}{3}M_{SH} (a,b) \leqslant \frac{1}{2}M_{AH}
(a,b)
 \leqslant \frac{1}{2}M_{SG} (a,b) \leqslant M_{AG} (a,b).
\end{equation}
\end{theorem}

\begin{proof} In order to prove the above theorem, we shall
prove each part separately.

Let us consider
\[
g_{SA\_SH} (x) = \frac{{f}''_{SA} (x)}{{f}''_{SH} (x)} = \frac{(x +
1)^3}{(x + 1)^3 + 4\sqrt 2 (x^2 + 1)^{3 / 2}}, \,\, x \in (0,\infty
),
\]

This gives
\begin{equation}
\label{eq39}
{g}'_{SA\_SH} (x) = - \frac{24(x - 1)(x^2 + 1)(x + 1)^2}{\sqrt {2(x^2 + 1)}
\left[ {(x + 1)^3 + 4\sqrt 2 (x^2 + 1)^{3 / 2}} \right]^2}
\begin{cases}
 { \geqslant 0,} & {x \leqslant 1} \\
 { \leqslant 0,} & {x \geqslant 1} \\
\end{cases}.
\end{equation}

In view of (\ref{eq39}) we conclude that the function $g_{SA\_SH}
(x)$ increasing in $x \in (0,1)$ and decreasing in $x \in (1,\infty
)$, and hence
\begin{equation}
\label{eq40}
\beta = \mathop {\sup }\limits_{x \in (0,\infty )} g_{SA\_SH} (x) =
g_{SA\_SH} (1) = \frac{1}{3}.
\end{equation}

Applying (\ref{eq33}) for the \textit{difference of means} $M_{SA}
(a,b)$ and $M_{SH} (a,b)$, and using (\ref{eq40}), we get
\begin{equation}
\label{eq41}
M_{SA} (a,b) \leqslant \frac{1}{3}M_{SH} (a,b).
\end{equation}

Let us consider
\[
g_{SH\_AH} (x) = \frac{{f}''_{SH} (x)}{{f}''_{AH} (x)} = \frac{(x +
1)^3 + 4\sqrt 2 (x^2 + 1)^{3 / 2}}{4\sqrt 2 (x^2 + 1)^{3 / 2}}, \,\,
x \in (0,\infty ),
\]

This gives
\begin{equation}
\label{eq42}
{g}'_{SH\_AH} (x) = - \frac{3(x - 1)(x + 1)^2}{4\sqrt 2 (x^2 + 1)^{5 / 2}}
\begin{cases}
 { \geqslant 0,} & {x \leqslant 1} \\
 { \leqslant 0,} & {x \geqslant 1} \\
\end{cases}.
\end{equation}

In view of (\ref{eq42}), we conclude that the function $g_{SH\_AH}
(x)$ is increasing in $x \in (0,1)$ and decreasing in $x \in
(1,\infty )$, and hence
\begin{equation}
\label{eq43}
\beta = \mathop {\sup }\limits_{x \in (0,\infty )} g_{SH\_AH} (x) =
g_{SG\_AH} (1) = \frac{3}{2}.
\end{equation}

Applying (\ref{eq33}) for the \textit{difference of means} $M_{SH}
(a,b)$ and $M_{AH} (a,b)$, and using (\ref{eq43}), we get
\begin{equation}
\label{eq44}
M_{SH} (a,b) \leqslant \frac{3}{2}M_{AH} (a,b).
\end{equation}

Let us consider
\[
g_{SG\_AH} (x) = \frac{{f}''_{SG} (x)}{{f}''_{AH} (x)} = \frac{(x +
1)^3\left[ {4x^{3 / 2} + \sqrt 2 \mbox{ }(x^2 + 1)^{3 / 2}}
\right]}{16\sqrt 2 \mbox{ }(x^2 + 1)^{3 / 2}x^{3 / 2}}, \,\, x \in
(0,\infty ),
\]

This gives
\begin{equation}
\label{eq45}
{g}'_{SG\_AH} (x) = \frac{3(x + 1)^4(x - 1)\left[ {\sqrt 2 \left( {x^2 + 1}
\right)^{5 / 2} - 8x^{5 / 2}} \right]}{32\sqrt 2 \left[ {x(x + 1)}
\right]^{5 / 2}}
\begin{cases}
 { \geqslant 0,} & {x \geqslant 1,} \\
 { \leqslant 0,} & {x \leqslant 1,} \\
\end{cases}
\end{equation}

\noindent where we have used the fact that $x^2 + 1 \geqslant 2x$,
$\forall x \in (0,\infty )$.

In view of (\ref{eq45}), we conclude that the function $g_{SG\_AH}
(x)$ is decreasing in $x \in (0,1)$ and increasing in $x \in
(1,\infty )$, and hence
\begin{equation}
\label{eq46}
\alpha = \mathop {\inf }\limits_{x \in (0,\infty )} g_{SG\_AH} (x) =
g_{SG\_AH} (1) = 1.
\end{equation}

Applying (\ref{eq33}) for the \textit{difference of means} $M_{AH}
(a,b)$ and $M_{SG} (a,b)$, and using (\ref{eq46}), we get
\begin{equation}
\label{eq47}
M_{AH} (a,b) \leqslant M_{SG} (a,b).
\end{equation}

Let us consider
\[
g_{SG\_AG} (x) = \frac{{f}''_{SG} (x)}{{f}''_{AG} (x)} = \frac{4x^{3
/ 2} + \sqrt 2 \mbox{ }(x^2 + 1)^{3 / 2}}{\sqrt 2 \mbox{ }(x^2 +
1)^{3 / 2}}, \,\, x \in (0,\infty ),
\]

This gives
\begin{equation}
\label{eq48}
{g}'_{SG\_AG} (x) = - \frac{6(x - 1)(x + 1)\sqrt x }{\sqrt 2 (x^2 + 1)^{5 /
2}}
\begin{cases}
 { \geqslant 0,} & {x \leqslant 1} \\
 { \leqslant 0,} & {x \geqslant 1} \\
\end{cases}.
\end{equation}

In view of (\ref{eq48}), we conclude that the function $g_{SG\_AG}
(x)$ is increasing in $x \in (0,1)$ and decreasing in $x \in
(1,\infty )$, and hence
\begin{equation}
\label{eq49}
M = \mathop {\sup }\limits_{x \in (0,\infty )} g_{SG\_AG} (x) = g_{SG\_AG}
(1) = 2.
\end{equation}

Applying (\ref{eq33}) for the \textit{difference of means} $M_{SG}
(a,b)$ and $M_{AG} (a,b)$, and using (\ref{eq46}), we get
\begin{equation}
\label{eq50}
\frac{1}{2}M_{SG} (a,b) \leqslant M_{AG} (a,b)
\end{equation}

Combining the results (\ref{eq41}), (\ref{eq44}), (\ref{eq47}) and
(\ref{eq50}) we get the proof of the inequality (\ref{eq38}).
\end{proof}

\begin{corollary} \label{cor31} The following inequalities hold:
\begin{align}
\label{eq51} & H(a,b) \leqslant G(a,b) \leqslant \frac{2H(a,b) +
S(a,b)}{3} \leqslant \frac{A(a,b) + H(a,b)}{2}\\
& \qquad \leqslant \frac{S(a,b) + G(a,b)}{2} \leqslant \frac{H(a,b)
+ 2S(a,b)}{3}\notag\\
& \qquad \qquad \leqslant A(a,b) \leqslant S(a,b) + H(a,b) -
G(a,b)\notag\\
& \qquad \qquad \qquad \leqslant S(a,b) \leqslant 3\left[ {A(a,b) -
G(a,b)} \right] +
 H(a,b).\notag
\end{align}
\end{corollary}

\begin{proof} Simplifying the results given in (\ref{eq41}),
(\ref{eq44}), (\ref{eq47}) and (\ref{eq50}) we get the required
result.
\end{proof}

\begin{remark} The inequalities (\ref{eq51}) are the
improvement over the following well known result
\begin{equation}
\label{eq52}
\min \left\{ {a,b} \right\} \leqslant H(a,b) \leqslant G(a,b) \leqslant
A(a,b) \leqslant S(a,b) \leqslant \max \left\{ {a,b} \right\}.
\end{equation}
\end{remark}

In the following corollary, we shall give a further improvement over
the inequalities (\ref{eq51}).

\begin{corollary} \label{cor32} The following inequalities hold:
\begin{align}
\label{eq53} & H(a,b) \leqslant \frac{2A(a,b)H(a,b)}{A(a,b) +
H(a,b)}
\leqslant G(a,b) \leqslant \frac{2H(a,b) + S(a,b)}{3}\\
& \quad \leqslant \frac{A(a,b) + H(a,b)}{2} \leqslant \sqrt
{\frac{\left( {A(a,b)} \right)^2 + \left( {H(a,b)} \right)^2}{2}}
\leqslant
\frac{S(a,b) + G(a,b)}{2}\notag\\
& \quad \quad \leqslant \frac{H(a,b) + 2S(a,b)}{3} \leqslant A(a,b)
\leqslant
S(a,b) + H(a,b) - G(a,b)\notag\\
& \quad \quad \quad \leqslant S(a,b) \leqslant 3\left[ {A(a,b) -
G(a,b)} \right] + H(a,b).\notag
\end{align}
\end{corollary}

\begin{proof} Replace $a$ by $A(a,b)$ and $b$ by $H(a,b)$ in
(\ref{eq52}) we get
\begin{align}
& \min \left\{ {A(a,b),H(a,b)} \right\} \leqslant H\left(
{A(a,b),H(a,b)} \right) \leqslant G\left( {A(a,b),H(a,b)}
\right)\notag\\
& \qquad \leqslant A\left( {A(a,b),H(a,b)} \right) \leqslant S\left(
 {A(a,b),H(a,b)}\right) \leqslant \max \left\{ {A(a,b),H(a,b)}
 \right\}.\notag
\end{align}

This gives
\begin{equation}
\label{eq54}
H(a,b) \leqslant \frac{2A(a,b)H(a,b)}{A(a,b) + H(a,b)} \leqslant G(a,b)
\leqslant \frac{A(a,b) + H(a,b)}{2}
\end{equation}
\[
 \leqslant \sqrt {\frac{\left( {A(a,b)} \right)^2 + \left( {H(a,b)}
\right)^2}{2}} \leqslant A(a,b) \leqslant S(a,b)
\]

The inequality (\ref{eq54}) gives a different kind of improvement
over the inequality (\ref{eq52}).

Let us consider
\begin{align}
\label{eq55} K(a,b) & = \frac{S(a,b) + G(a,b)}{2} - \sqrt
{\frac{A(a,b)^2 + H(a,b)^2}{2}}\\
& = \frac{\left( {\frac{S(a,b) + G(a,b)}{2}} \right)^2 -
\frac{A(a,b)^2 + H(a,b)^2}{2}}{\frac{S(a,b) + G(a,b)}{2} + \sqrt
{\frac{A(a,b)^2 + H(a,b)^2}{2}} }.\notag
\end{align}

Now we shall show that
\begin{equation}
\label{eq56}
\left( {\frac{S(a,b) + G(a,b)}{2}} \right)^2 - \frac{A(a,b)^2 + H(a,b)^2}{2}
\geqslant 0.
\end{equation}

For it, let us consider
\begin{align}
\label{eq57} k(x) & = \left[ {\frac{\sqrt {2(x^2 + 1)} }{4} +
\frac{\sqrt x }{2}} \right]^2 - \frac{1}{2}\left[ {\left( {\frac{x +
1}{2}} \right)^2 + \left( {\frac{2x}{x + 1}} \right)^2} \right]\\
& = \frac{8x^2 + (x + 1)^2\sqrt {2x(x^2 + 1)} }{4(x + 1)^2}
> 0, \,\, \forall x \in (0,\infty ).\notag
\end{align}

Now the expression (\ref{eq57}) together with (\ref{eq20}) give us
(\ref{eq56}), or equivalently, we can say that
\begin{equation}
\label{eq58}
\sqrt {\frac{A(a,b)^2 + H(a,b)^2}{2}} \leqslant \frac{S(a,b) + G(a,b)}{2}.
\end{equation}

Finally, the inequalities (\ref{eq51}), (\ref{eq54}) and
(\ref{eq58}) give us the proof of the inequalities (\ref{eq54}).
This completes the proof of the corollary.
\end{proof}

\begin{theorem} \label{the32} The following inequalities hold:
\begin{equation}
\label{eq59}
\frac{1}{8}M_{AH} (a,b) \leqslant M_{N_2 N_1 } (a,b) \leqslant
\frac{1}{3}M_{N_2 G} (a,b)
 \leqslant \frac{1}{4}M_{AG} (a,b) \leqslant M_{AN_2 } (a,b).
\end{equation}
\end{theorem}

\begin{proof} In order to prove the above theorem, we shall
prove each part separately.

Let us consider
\[
g_{AH\_N_2 N_1 } (x) = \frac{{f}''_{AH} (x)}{{f}''_{N_2 N_1 } (x)} =
\frac{32x^{5 / 2}(2x + 2)^{3 / 2}}{(x + 1)^3\left[ { - 2x - 2x^{5 /
2} + x(2x + 2)^{3 / 2}} \right]}, \,\, x \in (0,\infty ).
\]

This gives
\begin{align}
{g}'_{AH\_N_2 N_1 } (x) & = - \frac{48\sqrt {2x(x + 1)} }{(x +
1)^4\left[ { - 2x - 2x^{5 / 2} + x(2x + 2)^{3 / 2}}
\right]^2}\times\notag\\
& \qquad \times \left[ {4x^2(1 - x^{5 / 2}) + x^2(x - 1)(2x + 2)^{5
/ 2}} \right]\notag\\
& = \frac{48x^2(x + 1)\left( {1 - \sqrt x } \right)\sqrt {2x(x + 1)}
}{(x + 1)^4\left[ { - 2x - 2x^{5 / 2} + x(2x + 2)^{3 / 2}}
\right]^2}\times\notag\\
& \qquad \times \left[ {\sqrt 2 \left( {\sqrt x + 1} \right)\left(
{x + 1} \right)^{3 / 2} - \left( {x^2 + x^{3 / 2} + x + \sqrt x + 1}
\right)} \right].\notag
\end{align}

Since $\sqrt {2(x + 1)} \geqslant \sqrt x + 1$, $\forall x \in
(0,\infty )$, then this implies that
\begin{align}
\sqrt 2 (x + 1)^{3 / 2}\left( {\sqrt x + 1} \right) & \geqslant
\left( {\sqrt x + 1} \right)^2(x + 1)\notag\\
& \geqslant x^2 + x^{3 / 2} + x + \sqrt x + 1\notag
\end{align}

Thus we conclude that
\begin{equation}
\label{eq60} {g}'_{AH\_N_2 N_1 } (x)\begin{cases}
 { < 0,} & {x > 1} \\
 { > 0,} & {x < 1} \\
\end{cases}.
\end{equation}

In view of (\ref{eq60}), we conclude that the function $g_{AH\_N_2
N_1 } (x)$ is increasing in $x \in (0,1)$ and decreasing in $x \in
(1,\infty )$, and hence
\begin{equation}
\label{eq61}
\beta = \mathop {\sup }\limits_{x \in (0,\infty )} g_{AH\_N_2 N_1 } (x) =
g_{AH\_N_2 N_1 } (1) = 8.
\end{equation}

Applying (\ref{eq33}) for the \textit{difference of means} $M_{AH}
(a,b)$ and $M_{N_2 N_1 } (a,b)$ along with (\ref{eq61}), we get
\begin{equation}
\label{eq62}
\frac{1}{8}M_{AH} (a,b) \leqslant M_{N_2 N_1 } (a,b)
\end{equation}

Let us consider
\[
g_{N_2 N_1 \_N_2 G} (x) = \frac{{f}''_{N_2 N_1 } (x)}{{f}''_{N_2 G}
(x)} = \frac{ - 2x - 2x^{5 / 2} + x(2x + 2)^{3 / 2}}{2x\left[ {1 +
x^{3 / 2} - (2x + 2)^{3 / 2}} \right]}, \,\, x \in (0,\infty ).
\]

This gives
\begin{equation}
\label{eq63}
{g}'_{N_2 N_1 \_N_2 G_1 } (x) = \frac{3x^2\sqrt {2x + 2} \left( {1 - \sqrt x
} \right)}{2x^2\left[ { - 1 - x^{3 / 2} + (2x + 2)^{3 / 2}} \right]^2}
\begin{cases}
 { < 0,} & {x > 1,} \\
 { > 0,} & {x < 1.} \\
\end{cases}
\end{equation}

In view of (\ref{eq63}), we conclude that the function $g_{N_2 N_1
\_N_2 G} (x)$ is increasing in $x \in (0,1)$ and decreasing in $x
\in (1,\infty )$, and hence
\begin{equation}
\label{eq64}
\beta = \mathop {\sup }\limits_{x \in (0,\infty )} g_{N_2 N_1 \_N_2 G} (x) =
g_{N_2 N_1 \_N_2 G} (1) = \frac{1}{3}.
\end{equation}

Applying (\ref{eq33}) for the \textit{difference of means} $M_{N_2
N_1 } (a,b)$ and $M_{N_2 G} (a,b)$ along with (\ref{eq64}), we get
\begin{equation}
\label{eq65}
M_{N_2 N_1 } (a,b) \leqslant \frac{1}{3}M_{N_2 G} (a,b).
\end{equation}

Let us consider
\[
g_{N_2 G\_AG} (x) = \frac{{f}''_{N_2 G} (x)}{{f}''_{AG} (x)} = -
\frac{1 + x^{3 / 2} - (2x + 2)^{3 / 2}}{(2x + 2)^{3 / 2}}, \,\, x
\in (0,\infty ).
\]

This gives
\begin{equation}
\label{eq66}
{g}'_{N_2 G\_AG} (x) = \frac{3\left( {1 - \sqrt x } \right)}{(2x + 2)^{5 /
2}}
\begin{cases}
 { \leqslant 0,} & {x \geqslant 1} \\
 { \geqslant 0,} & {x \leqslant 1} \\
\end{cases}.
\end{equation}

In view of (\ref{eq66}), we conclude that the function $g_{AH\_N_2
N_1 } (x)$ is increasing in $x \in (0,1)$ and decreasing in $x \in
(1,\infty )$, and hence
\begin{equation}
\label{eq67}
\beta = \mathop {\sup }\limits_{x \in (0,\infty )} g_{N_2 G\_AG} (x) =
g_{N_2 G\_AG} (1) = \frac{3}{4}.
\end{equation}

Applying (\ref{eq33}) for the \textit{difference of means} $M_{N_2
G} (a,b)$ and $M_{AG} (a,b)$ along with (\ref{eq67}) we get
\begin{equation}
\label{eq68}
M_{N_2 G} (a,b) \leqslant \frac{3}{4}M_{AG} (a,b).
\end{equation}

Let us consider
\[
g_{AG\_AN_2 } (x) = \frac{{f}''_{AG} (x)}{{f}''_{AN_2 } (x)} =
\frac{(2x + 2)^{3 / 2}}{\left( {\sqrt x + 1} \right)\left( {x -
\sqrt x + 1} \right)}, \,\, x \in (0,\infty ).
\]

This gives
\begin{equation}
\label{eq69}
{g}'_{AG\_AN_2 } (x) = \frac{3\left( {1 - \sqrt x } \right)\sqrt {2x + 2}
}{\left( {\sqrt x + 1} \right)^2\left( {x - \sqrt x + 1} \right)^2}
\begin{cases}
 { \leqslant 0,} & {x \geqslant 1} \\
 { \geqslant 0,} & {x \leqslant 1} \\
\end{cases}.
\end{equation}

In view of (\ref{eq69}), we conclude that the function $g_{AG\_AN_2
} (x)$ is increasing in $x \in (0,1)$ and decreasing in $x \in
(1,\infty )$, and hence
\begin{equation}
\label{eq70}
\beta = \mathop {\sup }\limits_{x \in (0,\infty )} g_{AG\_AN_2 } (x) =
g_{AG\_AN_2 } (1) = 4.
\end{equation}

Applying (\ref{eq33}) for the \textit{difference of means} $M_{AG}
(a,b)$ and $M_{AN_2 } (a,b)$ along with (\ref{eq70}) we get the
required result.
\begin{equation}
\label{eq71}
\frac{1}{4}M_{AG} (a,b) \leqslant M_{AN_2 } (a,b).
\end{equation}

Combining the results (\ref{eq62}), (\ref{eq65}), (\ref{eq68}) and
(\ref{eq71}) we get the proof of the inequalities (\ref{eq59}).
\end{proof}

\begin{corollary} \label{cor33} The inequalities hold:
\begin{align}
\label{eq72} & H(a,b) \leqslant G(a,b) \leqslant \frac{G(a,b) +
H(a,b) + 3N_2 (a,b)}{5}\\
&\quad \leqslant \frac{G(a,b) + 2N_2 (a,b)}{3} \leqslant N_1 (a,b)
\leqslant \frac{2A(a,b) + 7N_1 (a,b)}{9} \leqslant N_2 (a,b)\notag\\
& \quad \quad \leqslant \frac{A(a,b) + N_1 (a,b)}{2} \leqslant
\frac{7A(a,b) + H(a,b)}{8} \leqslant A(a,b).\notag
\end{align}
\end{corollary}

\begin{proof} Follows in view of (\ref{eq57}), (\ref{eq60}),
(\ref{eq63}), (\ref{eq66}) and (\ref{eq31}).
\end{proof}

\begin{remark} \label{rem32} The inequalities (\ref{eq72}) can be
considered as an improvement over the following inequalities:
\begin{equation}
\label{eq73}
H(a,b) \leqslant G(a,b) \leqslant N_1 (a,b) \leqslant N_2 (a,b) \leqslant
A(a,b).
\end{equation}
\end{remark}

\begin{theorem} \label{the33} The following inequalities hold:
\begin{equation}
\label{eq74}
M_{SA} (a,b) \leqslant \frac{4}{5}M_{SN_2 } (a,b) \leqslant 4M_{AN_2 }
(a,b),
\end{equation}
\begin{equation}
\label{eq75}
M_{SH} (a,b) \leqslant 2M_{SN_1 } (a,b) \leqslant \frac{3}{2}M_{SG} (a,b),
\end{equation}

\noindent and
\begin{equation}
\label{eq76}
M_{SA} (a,b) \leqslant \frac{3}{4}M_{SN_3 } (a,b) \leqslant
\frac{2}{3}M_{SN_1 } (a,b).
\end{equation}
\end{theorem}

\begin{proof} In order to prove the above theorem, we shall
prove each part separately.

Let us consider
\[
g_{SA\_SN_2 } (x) = \frac{{f}''_{SA} (x)}{{f}''_{SN_2 } (x)}
 = \frac{8x^{3 / 2}(2x + 2)^{3 / 2}}{8x^{3 / 2}(2x + 2)^{3 / 2} + (1 + x^{3
/ 2})(2x^{3 / 2} + 2)^{3 / 2}}, \,\, x \in (0,\infty ).
\]

This gives
\begin{align}
{g}'_{SA\_SN_2 } (x) & = - \frac{96\sqrt {x(x^2 + 1)(x + 1)}
}{\left[ {8x^{3 / 2}(2x + 2)^{3 / 2} + (1 + x^{3 / 2})(2x^{3 / 2} +
2)^{3 / 2}} \right]^2}\times\notag\\
& \qquad \times \left[ {(x^2 - 1)(x^{5 / 2} + 1) + 2x(x^{5 / 2} -
1)} \right]\notag\\
& = - \frac{96\left( {\sqrt x - 1} \right)\sqrt {x(x^2 + 1)(x + 1)}
}{\left[ {8x^{3 / 2}(2x + 2)^{3 / 2} + (1 + x^{3 / 2})(2x^{3 / 2} +
2)^{3 / 2}} \right]^2}\times\notag\\
& \qquad \times \left[ {\left( {\sqrt x + 1} \right)(x + 1)(x^{5 /
2} + 1) + 2x\left( {x^2 + x^{3 / 2} + x + \sqrt x + 1} \right)}
\right].\notag
\end{align}

Thus, we have
\begin{equation}
\label{eq77} {g}'_{SA\_SN_2 } (x)\begin{cases}
 { \geqslant 0,} & {x \geqslant 1,} \\
 { \leqslant 0,} & {x \leqslant 1.} \\
\end{cases}
\end{equation}

In view of (\ref{eq77}), we conclude that the function $g_{SA\_SN_2
} (x)$ is increasing in $x \in (0,1)$ and decreasing in $x \in
(1,\infty )$, and hence
\begin{equation}
\label{eq78}
\beta = \mathop {\sup }\limits_{x \in (0,\infty )} g_{SA\_SN_2 } (x) =
g_{SA\_SN_2 } (1) = \frac{4}{5}.
\end{equation}

Applying (\ref{eq33}) for the \textit{difference of means} $M_{SA}
(a,b)$ and $M_{SN_2 } (a,b)$ along with (\ref{eq78}) we get
\begin{equation}
\label{eq79}
M_{SA} (a,b) \leqslant \frac{4}{5}M_{SN_2 } (a,b).
\end{equation}

Let us consider
\[
g_{SN_2 \_AN_2 } (x) = \frac{{f}''_{SN_2 } (x)}{{f}''_{AN_2 } (x)} =
\frac{8x^{3 / 2}(2x + 2)^{3 / 2} + (1 + x^{3 / 2})(2x^2 + 2)^{3 /
2}}{(2x^2 + 2)(x^{3 / 2} + 1)}, \,\, x \in (0,\infty ),
\]

This gives
\begin{align}
{g}'_{SN_2 \_AN_2 } (x) & = - \frac{12\left[ {x(x + 1)} \right]^{9 /
2}\left[ {(x^2 - 1)(1 + x^{5 / 2}) + 2x(x^{5 / 2} - 1)}
\right]}{(x^2 + 1)^{5 / 2}x^4(x + 1)^4(x^{3 / 2} + 1)^2}\notag\\
& = - \frac{12\left( {x(x + 1)} \right)^{9 / 2}\left( {\sqrt x - 1}
\right)}{(x^2 + 1)^{5 / 2}x^4(x + 1)^4(x^{3 / 2} +
1)^2}\times\notag\\
& \qquad \times \left[ {\left( {\sqrt x + 1} \right)(x + 1)(x^{5 /
2} + 1) + 2x\left( {x^2 + x^{3 / 2} + x + \sqrt x + 1} \right)}
\right].\notag
\end{align}

Thus we have
\begin{equation}
\label{eq80} {g}'_{SN_2 \_AN_2 } (x)\begin{cases}
 { \geqslant 0,} & {x \leqslant 1,} \\
 { \leqslant 0,} & {x \geqslant 1.} \\
\end{cases}
\end{equation}

In view of (\ref{eq80}), we conclude that the function $g_{SN_2
\_AN_2 } (x)$ is increasing in $x \in (0,1)$ and decreasing in $x
\in (1,\infty )$, and hence
\begin{equation}
\label{eq81}
\beta = \mathop {\sup }\limits_{x \in (0,\infty )} g_{SN_2 \_AN_2 } (x) =
g_{SN_2 \_AN_2 } (1) = \frac{4}{5}.
\end{equation}

Applying (\ref{eq33}) for the \textit{difference of means} $M_{SN_2
} (a,b)$ and $M_{AN_2 } (a,b)$ along with (\ref{eq81}) we get
\begin{equation}
\label{eq82}
\frac{1}{5}M_{SN_2 } (a,b) \leqslant M_{AN_2 } (a,b).
\end{equation}

Combining the results (\ref{eq79}) and (\ref{eq82}) we get the proof
of the inequalities (\ref{eq74}). Now we shall give the proof of
(\ref{eq75}).

Let us consider
\[
g_{SH\_SN_1 } (x) = \frac{{f}''_{SH} (x)}{{f}''_{SN_1 } (x)} =
\frac{16x^{3 / 2}\left[ {(x + 1)^3 + 2(2x^2 + 2)^{3 / 2}}
\right]}{(x + 1)^3\left[ {16x^{3 / 2} + (2x^2 + 2)^{3 / 2}}
\right]}, \,\, x \in (0,\infty ),
\]

This gives
\begin{align}
{g}'_{SH\_SN_1 } (x) & = - \frac{48\sqrt {2x^2 + 2} }{x^2(x +
1)^4\left[ {16x^{3 / 2} + (2x^2 + 2)^{3 / 2}}
\right]^2}\times\notag\\
& \qquad \times \left[ {64x^{9 / 2}(1 - x) + 5x^4(x^2 - 1) +
4x^3(x^4 - 1)} \right.\notag\\
& \qquad \qquad \left. { + x^2(x^2 - 1) + x^2(x - 1)(2x^2 + 2)^{5 /
2}}
\right]\notag\\
& = - \frac{1536x^2(x - 1)\sqrt {2x^2 + 2} }{x^2(x + 1)^4\left[
{16x^{3 / 2} + (2x^2 + 2)^{3 / 2}} \right]^2}\times\notag\\
& \qquad \times \left\{ {\left[ {\left( {\frac{x + 1}{2}} \right)^5
- \left( {\sqrt x } \right)^5} \right] + \left[ {\left( {\sqrt
{\frac{x^2 + 1}{2}} } \right)^5 - \left( {\sqrt x } \right)^5}
\right]} \right\}.\notag
\end{align}

Since $S(a,b) \geqslant A(a,b) \geqslant G(a,b)$, one gets
\begin{equation}
\label{eq83} {g}'_{SH\_SN_1 } (x)\begin{cases}
 { \geqslant 0,} & {x \leqslant 1,} \\
 { \leqslant 0,} & {x \geqslant 1.} \\
\end{cases}
\end{equation}

In view of (\ref{eq83}) we conclude that the function $g_{SH\_SN_1 }
(x)$ increasing in $x \in (0,1)$ and decreasing in $x \in (1,\infty
)$, and hence
\begin{equation}
\label{eq84}
\beta = \mathop {\sup }\limits_{x \in (0,\infty )} g_{SH\_SN_1 } (x) =
g_{SH\_SN_1 } (1) = 2.
\end{equation}

Applying (\ref{eq33}) for the \textit{difference of means} $M_{SH}
(a,b)$ and $M_{SN_1 } (a,b)$ along with (\ref{eq84}) we get
\begin{equation}
\label{eq85}
M_{SH} (a,b) \leqslant 2M_{SN_1 } (a,b).
\end{equation}

Let us consider
\[
g_{SN_1 \_SG} (x) = \frac{{f}''_{SN_1 } (x)}{{f}''_{SG} (x)} =
\frac{8x^{3 / 2} + (x^2 + 1)\sqrt {2x^2 + 2} }{2\left[ {4x^{3 / 2} +
(x^2 + 1)\sqrt {2x^2 + 2} } \right]}, \,\, x \in (0,\infty ),
\]

This gives
\begin{equation}
\label{eq86}
{g}'_{SN_1 \_SG} (x) = - \frac{3(x^2 - 1)\sqrt {2x^3 + 2x} }{\left[ {4x^{3 /
2} + (x^2 + 1)\sqrt {2x^2 + 2} } \right]^2}
\begin{cases}
 { \geqslant 0,} & {x \leqslant 1,} \\
 { \leqslant 0,} & {x \geqslant 1.} \\
\end{cases}
\end{equation}

In view of (\ref{eq86}), we conclude that the function $g_{SN_1
\_SG} (x)$ is increasing in $x \in (0,1)$ and decreasing in $x \in
(1,\infty )$, and hence
\begin{equation}
\label{eq87}
\beta = \mathop {\sup }\limits_{x \in (0,\infty )} g_{SN_1 \_SG} (x) =
g_{SN_1 \_SG} (1) = \frac{3}{4}.
\end{equation}

Applying (\ref{eq33}) for the \textit{difference of means} $M_{SN_1
} (a,b)$ and $M_{SG} (a,b)$ along with (\ref{eq87}) we get
\begin{equation}
\label{eq88}
M_{SN_1 } (a,b) \leqslant \frac{3}{4}M_{SG} (a,b).
\end{equation}

Combining the results given in (\ref{eq85}) and (\ref{eq88}) we get
the proof of the inequalities (\ref{eq75}). Let us prove now the
inequalities (\ref{eq76}).

Let us consider
\[
g_{SA\_SN_3 } (x) = \frac{{f}''_{SA} (x)}{{f}''_{SN_3 } (x)} =
\frac{24x^{3 / 2}}{24x^{3 / 2} + (2x^2 + 2)^{3 / 2}}, \,\, x \in
(0,\infty ),
\]

This gives
\begin{equation}
\label{eq89}
{g}'_{SA\_SN_3 } (x) = - \frac{72(x - 1)(x + 1)\sqrt {2x(x^2 + 1)} }{\left[
{24x^{3 / 2} + (2x^2 + 2)^{3 / 2}} \right]^2}
\begin{cases}
 { \geqslant 0,} & {x \leqslant 1,} \\
 { \leqslant 0,} & {x \geqslant 1.} \\
\end{cases}
\end{equation}

In view of (\ref{eq89}), we conclude that the function $g_{SA\_SN_3
} (x)$ is increasing in $x \in (0,1)$ and decreasing in $x \in
(1,\infty )$, and hence
\begin{equation}
\label{eq90}
\beta = \mathop {\sup }\limits_{x \in (0,\infty )} g_{SA\_SN_3 } (x) =
g_{SA\_SN_3 } (1) = \frac{3}{4}.
\end{equation}

Applying (\ref{eq33}) for the \textit{difference of means} $M_{SA}
(a,b)$ and $M_{SN_3 } (a,b)$ along with (\ref{eq90}) we get
\begin{equation}
\label{eq91}
M_{SA} (a,b) \leqslant \frac{3}{4}M_{SN_3 } (a,b).
\end{equation}

Let us consider
\[
g_{SN_3 \_SN_1 } (x) = \frac{{f}''_{SN_3 } (x)}{{f}''_{SN_1 } (x)} =
\frac{2\left[ {24x^{3 / 2} + (2x^2 + 2)^{3 / 2}} \right]}{3\left[
{16x^{3 / 2} + (2x^2 + 2)^{3 / 2}} \right]}, \,\, x \in (0,\infty ),
\]

This gives
\begin{equation}
\label{eq92}
{g}'_{SN_3 \_SN_1 } (x) = - \frac{(x - 1)(x + 1)\sqrt {2x(x^2 + 1)} }{\left[
{16x^{3 / 2} + (2x^2 + 2)^{3 / 2}} \right]^2}
\begin{cases}
 { \geqslant 0,} & {x \leqslant 1,} \\
 { \leqslant 0,} & {x \geqslant 1.} \\
\end{cases}
\end{equation}

In view of (\ref{eq92}), we conclude that the function $g_{SN_3
\_SN_1 } (x)$ is increasing in $x \in (0,1)$ and decreasing in $x
\in (1,\infty )$, and hence
\begin{equation}
\label{eq93}
\beta = \mathop {\sup }\limits_{x \in (0,\infty )} g_{SN_3 \_SN_1 } (x) =
g_{SN_3 \_SN_1 } (1) = \frac{3}{4}.
\end{equation}

Applying (\ref{eq33}) for the \textit{difference of means} $M_{SN_3
} (a,b)$ and $M_{SN_1 } (a,b)$ along with (\ref{eq93}) we get
\begin{equation}
\label{eq94}
M_{SN_3 } (a,b) \leqslant \frac{8}{9}M_{SN_1 } (a,b).
\end{equation}

Combining the results given in (\ref{eq91}) and (\ref{eq94}) we get
the proof of the inequalities (\ref{eq79}). This completes the proof
of the theorem.
\end{proof}

\begin{corollary} \label{cor34} The following inequalities hold:
\begin{align}
\label{eq95} & G(a,b) \leqslant \frac{S(a,b) + 3G(a,b)}{4} \leqslant
N_1 (a,b) \leqslant \frac{S(a,b) + 8N_1 (a,b)}{9}\\
&\quad \leqslant N_3 (a,b) \leqslant N_2 (a,b) \leqslant
\frac{A(a,b) + N_1 (a,b)}{2} \leqslant \frac{S(a,b) + 2N_1
(a,b)}{3}\notag\\
& \quad \quad \leqslant \left( {\frac{S(a,b) + 4N_2 (a,b)}{5}}
\right)\;or\;\left( {\frac{S(a,b) + 3N_3 (a,b)}{4}} \right)
\leqslant A(a,b)\notag
\end{align}
\noindent and

\begin{equation}
\label{eq96} G(a,b) \leqslant \;\frac{S(a,b) + 2H(a,b)}{2} \leqslant
N_1 (a,b) \leqslant \;\frac{S(a,b) + H(a,b)}{2} \leqslant N_2 (a,b).
\end{equation}

\end{corollary}

\begin{proof} The inequalities (\ref{eq74})-(\ref{eq76}) lead
us to (\ref{eq95}) and (\ref{eq96}).
\end{proof}

\begin{remark} \label{rem33} The inequalities (\ref{eq95}) can be
considered as refinement over the inequality (\ref{eq7}). Thus we
have three different kind of refinements given by (\ref{eq53}),
(\ref{eq72}) and (\ref{eq95}) for the inequality (\ref{eq7}). The
inequalities (\ref{eq96}) gives alternative improvement among the
\textit{means} $G(a,b)$, $N_1 (a,b)$ and $N_2 (a,b)$.
\end{remark}


\begin{thebibliography}{99}
\setlength{\itemsep}{5pt}

\bibitem{beb} E.F. BECKENBACH and R. BELLMAN, Inequalities, Springer-Verlag,
New York, 1971.

\bibitem{bob} J. BORWEIN and P. BORWEIN, \textit{$\pi$ and the AGM},
Wiley, 1987.

\bibitem{drp} S.S. DRAGOMIR  and   C.E.M. PEARCE, Selected Topics on
Hermite-Hadamard Inequalities and Applications, Research Report
Collection, Monograph, 2002, available on line:
rgmia.vu.edu.au/monographs/index.html

\bibitem{tan1} I.J. TANEJA, On a Difference of Jensen Inequality
and its Applications to Mean Divergence Measures -- \textit{RGMIA
Research Report Collection}, \textit{http://rgmia.vu.edu.au},
\textbf{7}(4)(2004), Art. 16. Also in:arXiv:math.PR/0501302 v1 19
Jan 2005.

\bibitem{tan2} I.J. TANEJA, On Symmetric and Non-Symmetric Divergence
Measures and Their Generalizations, to appear as a chapter in:
Advances in Imaging and Electron Physics, 2005.

\bibitem{zhw} ZHI-HUA ZHANG and YU-DONG WU, The New Bounds of the
Logarithmic Mean, \textit{RGMIA Research Report Collection},
\textit{http://rgmia.vu.edu.au}, \textbf{7}(2)(2004), Art. 7.

\end{thebibliography}
\end{document}